\documentclass[10pt,twoside]{article}
\usepackage{graphicx}
\usepackage{amsmath}
\usepackage{amsthm}
\usepackage{amsfonts}
\usepackage{Latex-document}

\title{\bf  Conformal Invariance, Universality,
and \vskip -2mm the Dimension  of the Brownian Frontier\vskip 6mm}
\def \R {\mathbb R}
\def \C {\mathbb C}
\def \hdim {{\rm dim}_h}
\def \Prob{{\bf P}}
\def \H {\mathbb H}
\def \p {\partial}
\def \Disk {\mathbb D}
\def \E {{\bf E}}
\def \Z {\mathbb Z}
\author{G. Lawler\vspace*{-0.5cm}\thanks{Department of Mathematics,
Cornell University, Malott Hall, Ithaca, NY 14853-4201 and
Department of Mathematics, Duke University, Durham, NC 27708-0320,
USA. E-mail: lawler@math.cornell.edu}}
\date{\vspace{-8mm}}

\begin{document}

\maketitle

\thispagestyle{first} \setcounter{page}{63}

\begin{abstract}\vskip 3mm
This paper describes joint work with Oded Schramm and Wendelin
Werner establishing the values of the planar Brownian intersection
exponents from which one derives the
Hausdorff dimension of certain exceptional
sets of planar Brownian motion.  In particular, we proof a conjecture of
Mandelbrot that the dimension of the  frontier is $4/3$.
The proof uses a universality principle for conformally
invariant  measures and a new process, the
stochastic Loewner evolution ($SLE$), introduced by
Schramm.  These ideas can be used
to study other planar lattice models from statistical physics
at criticality.   I discuss applications
to critical percolation on the triangular lattice,   loop-erased
random walk, and self-avoiding walk.

\vskip 4.5mm


\end{abstract}

\vskip 12mm

\section{Exceptional sets for planar Brownian motion} \label{section 1}\setzero
\vskip-5mm \hspace{5mm }

Let $B_t$ be a standard Brownian motion taking values
in $\R^2 = \C$ and let $B[s,t]$ denote the random set
$ B[s,t] = \{B_r: s \leq r \leq t\} .$
For $0 \leq t \leq 1$, we say that $B_t$ is a
\begin{itemize}
\item {\em cut point} for $B[0,1]$ if $B[0,t) \cap B(t,1] = \emptyset $;
\item {\em frontier point} for $B[0,1]$ if $B_t$ is on the boundary of
  the unbounded component of $\C \setminus B[0,1]$;
\item {\em pioneer point} for $B[0,1]$ if $B_t$ is on the boundary of
  the unbounded component of $\C \setminus B[0,t]$,
   i.e., if $B_t$ is a frontier point for $B[0,t] $.
\end{itemize}
I will discuss the following result proved
by Oded Schramm, Wendelin Werner, and myself.

\medskip

{\bf Theorem 1.}    \cite{LSW1,LSW2,LSWanal}
\it   If $B_t$ is a standard Brownian motion
in $\R^2 = \C$, then with probability one,
\[  \hdim\; [\; \mbox{cut points for } B[0,1]\; ] = 3/4 , \]
\[   \hdim\; [\; \mbox{frontier points for } B[0,1]\; ]= 4/3 , \]
\[  \hdim \; [\; \mbox{pioneer points for } B[0,1]\;]= 7/4 , \]
where $\hdim$ denotes Hausdorff dimension.
\rm

\medskip

Mandelbrot \cite{mandel} first gave the conjecture for the
Brownian frontier, basing his conjecture on numerical
simulation and then noting that simulations of the frontier
resembled simulations of self-avoiding walks.  It is
conjectured that the scaling limit of planar self-avoiding walks
has paths of dimension $4/3$.  Duplantier
and Kwon \cite{DK1} used nonrigorous conformal field
theory techniques to make the above conjectures for
the cut points and pioneer points.  More precisely, they
made conjectures about certain exponents called
the Brownian or simple
random walk {\em intersection exponents}.  More recently,
Duplaniter \cite{D3} has given other nonrigorous
arguments for the conjectures using quantum gravity.

To prove Theorem 1,  it suffices  to
find the values of
the Brownian   intersection exponents.
In fact, before Theorem 1 had been proved, it had been
established \cite{Lcut,Lfront,LBuda} that the Hausdorff
dimensions of the set of cut points, frontier points, and
pioneer points were $2 - \eta_1, 2 - \eta_2,$ and
$2 - \eta_3$, respectively, where $\eta_1,\eta_2,\eta_3$
are defined by saying that as $\epsilon \rightarrow
0+$,
\[  \Prob\{B[0,\frac{1}{2} - \epsilon^2] \cap
   B[\frac12 + \epsilon^2 ,1] = \emptyset \} \approx
    \epsilon^{\eta_1},  \]
\[  \Prob\{B[0,\frac{1}{2} - \epsilon^2] \cup
       B[\frac{1}{2} + \epsilon^2,1] \mbox{ does not disconnect }
  B_{1/2} \mbox{ from infinity} \}  \approx
   \epsilon^{\eta_2},   \]
\[ \Prob\{B[\epsilon^2,1] \mbox{ does not disconnect }
  0 \mbox{ from infinity} \}  \approx
   \epsilon^{\eta_3}.\]
It had also been established \cite{BL2,LP} that the analogous
exponents for simple random walk are the same as for
Brownian motion.

There are two main ideas in the proof.  The first is a
one parameter family of conformally invariant processes
developed by Oded Schramm \cite{S2} which he named the
Stochastic Loewner evolution ($SLE$).  The second is the
idea of  ``universality''  which states roughly that all conformally
invariant measures  that satisfy a certain
``locality'' or
``restriction'' property must have the same exponents as
Brownian motion (see \cite{LW2}).
In this paper, I will define $SLE$ and give
some of its properties; describe how analysis of $SLE$ leads
to finding the Brownian intersection exponents; and  finally
describe some other planar lattice models in statistical physics
at criticality that can be
understood using $SLE$.

\section{Stochastic Loewner evolution} \label{section 2}
\setzero\vskip-5mm \hspace{5mm }

  I will give a brief introduction to the stochastic
Loewner evolution ($SLE$); for more details, see
\cite{S1,LSW1,LSW2, LVienna,SR}.
Let $W_t$ denote a standard one dimensional Brownian motion.
If $\kappa \geq 0$ and $z$ is in the upper half plane
$\H = \{w \in \C: \Im(w) > 0\}$, let $g_t(z)$ be the solution to
the Loewner differential equation
\begin{equation}  \label{lde.chord}
                 \p_t g_t(z) = \frac{2}{g_t(z) - \sqrt{\kappa} \;
              W_t }, \;\;\;\;  g_0(z) = z .
\end{equation}
For each $z \in \H$, the solution $g_t(z)$ is defined up to
a time $T_z \in (0,\infty]$.   Let  $H_t=
\{z : T_z > t\}$.  Then $g_t$ is the unique
conformal transformation of $H_t$ onto $\H$ with
$g_t(z) - z = o(1)$ as $z \rightarrow \infty$.  In fact,
\[           g_t(z) = z + \frac{2t}{z} + O(\frac{1}{|z|^2}),
 \;\;\;\; z \rightarrow \infty . \]
It is easy to show that the maps $g_t$ are well defined.
It has been shown \cite{SR,LSWlerw} that there is a (random)
continuous
path $\gamma:[0,\infty) \rightarrow \overline \H$ such that $H_t$
is the unbounded component of $\H \setminus \gamma[0,t]$ and
$g_t(\gamma(t)) = \sqrt \kappa \: W_t$.  The conformal maps
$g_t$ or the corresponding paths $\gamma(t)$ are called the
{\em  chordal  stochastic Loewner evolution with parameter
$\kappa$ (chordal $SLE_\kappa$)}.  It is easy to check that
the distribution of
$SLE_\kappa$ is invariant (modulo time change) under
dilations $z \mapsto rz$.  Using this, we can use conformal
transformations to  define
chordal $SLE_\kappa$ connecting two distinct boundary
points of any simply connected domain.  This gives a family
of probability measures on curves (modulo
reparametrization) on such domains that is invariant under
conformal transformation.

Chordal $SLE_\kappa$ can also be considered as the
only probability distributions on continuous curves (modulo
reparametrization) $\gamma:[0,\infty)  \rightarrow
\overline \H$ with the following properties.
\begin{itemize}
\item $\gamma(0) = 0, \gamma(t) \rightarrow \infty$
as $t \rightarrow \infty$, and
$\gamma(t) \in \p H_t$ for all $t \in [0,\infty),$ where
$H_t$ is the unbounded component of $\H \setminus
\gamma[0,t].$
\item Let $h_t: H_t \rightarrow \H$ be the  unique
conformal transformation  with
$h_t(\gamma(t)) = 0, h_t(\infty) = \infty, h_t'(\infty) = 1$. Then
the conditional distribution of $\hat \gamma(s)
:= h_t \circ \gamma(s+t), 0 \leq s < \infty$,
given $\gamma[0,t]$ is  the same
as the original distribution.
\item The measure is invariant under $x + iy \mapsto
 -x+iy .$
\end{itemize}

There is a similar process called {\em radial $SLE_\kappa$}
on the unit disk.  Let $W_t$ be as above, and for $z$
in the unit disk $\Disk$, consider the equation
\[          \p_t g_t(z) = g_t(z) \; \frac{e^{i \sqrt \kappa W_t} +
           g_t(z)}{e^{i \sqrt \kappa W_t} - g_t(z) } , \;\;\;\;
           g_0(z) = z . \]
Let $U_t$ be the set of $z \in \Disk$ for which   $g_t(z)$
is defined.
It can be shown that there is a random path $\gamma:[0,\infty)
\rightarrow \overline \Disk$, such that $U_t$ is the component of
$\Disk \setminus \gamma[0,t]$ containing the origin;
$g_t(\gamma(t)) = e^{i \sqrt \kappa W_t} $; and $g_t$ is
a conformal transformation of $U_t$ onto $\Disk$ with
$g_t(0) = 0, g_t'(0) = e^{t}$.  We can define radial
$SLE_\kappa$ connecting any boundary point to
any interior point of a simply connected domain by
conformal transformation.

The qualitative behavior of the paths $\gamma$ varies
considerably as $\kappa$ varies, although chordal
and radial $SLE_\kappa$ for the same $\kappa$
are qualitatively similar.  The
Hausdorff dimension of $\gamma[0,t]$ for chordal
or radial $SLE_{\kappa}$ is conjectured to be
$\min\{1 + (\kappa/8),2\}$.   This has been proved
for $\kappa = 8/3, 6,$ see \cite{B}, and for other
$\kappa$ it is a rigorous upper bound \cite{SR}.
For $0 \leq \kappa \leq 4$, the paths $\gamma$
are simple (no self-intersections) and $\gamma(0,\infty)$
is a subset of $\H$ or $\Disk$.  For $\kappa >4$, the
paths have double points and hit $\p \H$ or $\p \Disk$
infinitely often.  If $\kappa \geq 8$, the paths are
space filling.

Investigation  of $SLE_\kappa$ requires studying
the    behavior of $SLE_\kappa$ under  conformal
maps.
Suppose $A$ is a compact subset of $\overline
\H$ not containing the
origin such that $\overline{A \cap \H} = A$
and $\H \setminus A$ is simply connected.
Let $\Phi$ denote the conformal transformation of $\H
\setminus A$ onto $\H$ with $\Phi(0) = 0 , \Phi(\infty)
= \infty, \Phi'(\infty) = 1$.  Let $\gamma$ denote a chordal
$SLE_\kappa$ starting at the origin,
and let  $T$ be the first time $t$ that
$A \cap \H \not\subset H_t$.  For $t < T$, let
$\tilde \gamma(t)
= \Phi \circ \gamma(t)$.  Let $\tilde g_t$ be the conformal
transformation of the unbounded component of $\H
\setminus \tilde \gamma[0,t]$ onto $\H$ with $\tilde g_t(z) - z
=o(1)$ as $z \rightarrow \infty$; define $a(t)$ by
$\tilde g(z) - z \sim a(t) \; z^{-1} . $  Then $\tilde g(t)$
satisfies the modified Loewner equation
\[
   \p_t \tilde g(t) = \frac{ \p_t a} {\tilde g_t(z) - \tilde W_t} , \;\;\;\;
                 \tilde g_0(z) = z,
\]
for some $\tilde W_t$.   In fact $\tilde W_t =
\tilde g_t \circ \Phi \circ g_t^{-1} (\sqrt \kappa \: W_t)$.
Using  the Loewner differential
equation and It\^o's formula,  we can write $\tilde W_t$ as
a local semimartingale, $d \tilde W_t = b(t) \; dt +
\sqrt{ \kappa \; \p_t a/2 } \;
dW_t$; here $b(t)$ and $a(t)$ are random depending
on $W_s, 0 \leq s \leq t$.   For $\kappa = 6$,
and only $\kappa =6$,
the drift term $b(t)$ disappears and hence $\tilde W_t$ is
a time change of Brownian motion.

\medskip

{ \bf Locality property for {\boldmath $SLE_6$}.}  \cite{LSW1} \it
If $\kappa = 6$, $\tilde \gamma(t), 0 \leq t < T$, has the same
distribution as a time change of $SLE_6$. \rm

\medskip

For other values of $\kappa$, the image $\tilde \gamma(t),
t < T $, has
a distribution that is absolutely continuous with respect to
that of (a time change of)
$SLE_\kappa$.  This follows from Girsanov's theorem
(see, e..g, \cite[Theorem I.6.4]{Bass}) that states roughly that
Brownian motions with the same variance but different drifts
give rise to absolutely continuous measures on paths.
Similarly, radial $SLE_\kappa$ can be obtained from
chordal $SLE_\kappa$ by considering its image
under a map taking $\H$ to $\Disk$.
  For all values of
$\kappa$ we get absolutely continuous measures (which
is why radial $SLE_\kappa$ is qualitatively the same
as chordal $SLE_\kappa$), but for $\kappa = 6$ we get
a special relationship \cite[Theorem 4.1]{LSW2}.

One of the reasons that $SLE_\kappa$ is useful is
that ``crossing probabilities'' and ``critical exponents'' for
the process can be calculated.  The basic idea is
to relate an event about the planar path $\gamma$
to an event about the driving process $\sqrt{\kappa}
W_t$  and then to use standard methods of stochastic
calculus to relate this to solutions of partial
differential equations.  As an example, consider chordal
$SLE_6$ in the upper half plane $\H$ going from
$x \in (0,1)$ to infinity.  Let $T$ be the first time
$t$ that $\gamma(t) \in (-\infty,0] \cap [1,\infty)$;
since $\kappa > 4$, $T < \infty$ with probability one.
Let
${\cal E}  $
be the event that $\gamma(T) \in (-\infty,0]$,
 and let
$ H_T$ be the unbounded component of
$ \H \setminus \gamma[0,T]$.   Let
$y_1,y_2$ be the minimum and maximum
of $\gamma[0,T] \cap \R$; on the event
${\cal E}$, $y_1 \leq 0$ and $x \leq y_2 < 1$.
  Let ${\cal L}$ denote the
{\em $\pi$-extremal distance } between $(-\infty,y_1]$ and
$[y_2,1]$
in  $H_T$,
 i.e., the number ${\cal L}$
such that $H_T$ can be mapped conformally onto $[0,{\cal L}] \times
[0,\pi]$ in a way that $(-\infty,y_1]$ and $[y_2,1]  $   are mapped
 onto the vertical boundaries.
   In order to relate
$SLE_6$ to  intersection exponents for Brownian motion
one needs to understand  the behavior of
 $\E^x[1_{\cal E} \; \exp\{-\lambda {\cal L}\}]$
as $x \rightarrow 1-$ for $\lambda \geq 0$.  It
is not hard to show that this  quantity
is closely related to $\E^x[1_{\cal E} \; g_T'(1)^\lambda]$.    If
we differentiate (\ref{lde.chord}) with respect to $z$ we get
an equation for $\p_t g'_t(1)$,  and standard techniques
of stochastic calculus can be applied to give a differential
equation for the function $r(x, \lambda)
= \E^x[1_{\cal E} \; g_T'(1)^\lambda]$.
We get an exact solution in terms of hypergeometric
functions \cite[Theorem 3.2]{LSW1}.  If $\lambda =0$,
so that $r(x) = \Prob^x[{\cal E}]$, we get
the formula given
by Cardy \cite{Cardy} for crossing probabilities
of percolation clusters (see \S\ref{percsec}).

\section{Applications}

\vskip-5mm \hspace{5mm }

\subsection{Brownian motion}\vskip-5mm \hspace{5mm }

As already mentioned, computation of dimensions for many
exceptional sets for Brownian motion reduces to finding the
Brownian intersection exponents.     These
exponents, which can be defined in
terms of crossing probabilities for non-intersecting
paths,  were studied in \cite{LW1,LW2}.
In these papers,  relations
were given between different exponents and a ``universality''
principle was shown for conformally invariant processes satisfying
an additional hypothesis (the term completely conformally
invariant was used there).
 Heuristic arguments indicated that
self-avoiding walks and percolation should also
 satisfy this hypothesis.   Unfortunately, from a rigorous
standpoint, we had only reduced a hard problem, computing
the Brownian intersection exponents, to the even harder problem
of showing conformal invariance and computing the exponents
for self-avoiding walks or critical percolation.

At the same time Schramm \cite{S1} was completing his beautiful
 construction of  $SLE_\kappa$ and conjecturing that
$SLE_6$ gave the boundaries of critical percolation clusters.  While
he was unable to prove that critical percolation has a conformally
invariant limit, he was able to conclude that if the limit
was conformally invariant then it must be $SLE_6$.  The
identification $\kappa =6$ was determined from rigorous
``crossing probabilities'' for $SLE_\kappa$; only $\kappa = 6$
was consistent with Cardy's formula (see \S\ref{percsec}) or even
the simple fact that a square should have crossing probability
$1/2$.

Since both Brownian motion and $SLE_6$ were conjectured
to be related to the scaling limit of critical percolation, it
was natural to try to use $SLE_6$ to prove results
about Brownian motion
 (and, as mentioned
before, the Hausdorff dimension of exceptional sets on
the path); see \cite{LSW1,LSW2, LSWanal,LSW3}.
    There were two major parts of the proof.  First,
the locality property for $SLE_6$ was formulated and
proved; this allowed ideas as in \cite{LW2}
to show  that the exponents of $SLE_6$  can be used
to find the exponents for Brownian motion.
Second, the exponents for $SLE_6$ had to be computed.
The basic idea is discussed at the end of the last section.
What makes $SLE$ so powerful is that it reduces problems
about a two-dimensional process to analysis of a
one-dimensional stochastic differential equation (and
hence a partial differential equation in one space
variable).

The universality in these papers was in terms of exponents.
We now know that the paths of planar Brownian
motion and $SLE_6$ are even more closely related.  The
``hull'' generated by an $SLE_6$ is the same as
the hull generated by a Brownian motion with oblique
reflection (see \cite{W}).  In particular, the frontiers (outer
boundaries) of the two processes have the same dimension.
There are now direct proofs that the Hausdorff dimension
of the frontier of $SLE_6$ is $4/3$ (\cite{B}) and this
stronger universality principle implies the same holds for
Brownian paths.

\subsection{Critical percolation}  \label{percsec}\vskip-5mm \hspace{5mm }

Suppose each vertex
 of the planar triangular lattice is colored independently
white or black, with the
probability of a white being $1/2$.   This
is  called {\em critical percolation} (on the
triangular lattice).  Let $D$ be a simply
connected domain in $\C = \R^2$ and  let $A_1,
A_2$ be disjoint nontrivial
connected arcs on $\p D$.
Consider the limit as $\delta \rightarrow 0$
of the probability that in critical percolation
on a lattice with mesh size $\delta$ that there is
a connected set of white vertices in $D$
connecting $A_1,A_2$.  It has long been
believed that this limit, $p(A_1,A_2;D)$,
 exists and is strictly between $0$
and $1$.   (Note: if the probability of a
white vertex is $p$, then $p(A_1,A_2;D)$ is
$0$ for
$p < 1/2$ and $1$ for $p > 1/2$.  One of the
features of {\em critical} percolation is the
fact that this quantity is strictly between $0$ and $1$.)
Moreover, it has
been conjectured that $p(A_1,A_2;D)$
is a conformal invariant
\cite{Cardy,Lang} .   It is also believed that
this limit does not depend on the nature of the
lattice; for example,  critical bond percolation in
$\Z^2$ (each bond is colored white or black
independently with probability $1/2$) should give
the same limit.

Cardy \cite{Cardy} used nonrigorous methods from
conformal field theory to find an exact formula
for $p(A_1,A_2;D)$; his calculations were done
for $D = \H$ and the formula involves hypergeometric
functions. Carleson noted that the formula was
much nicer if one chooses $D$ to be an equilateral
triangle of side length $1$; $A_1$, one of the sides; and
$A_2$, a line segment of length $x$ with one
endpoint on the vertex opposite $A_1$.  In this
case, Cardy's formula is $p(A_1,A_2;D)
=x$.  Schramm \cite{S1} went further and,
 {\em assuming existence and
conformal invariance of the limit}, showed
that the limiting boundary between black and
white clusters can be given in terms of $SLE_6$.
If $A_3$ denotes the third side
of the triangle (so that $A_3 \cap A_2$ is a
single point), we can consider the limiting cluster
formed by taking all the white vertices that are
connected by a path of white vertices to $A_3$.
In the limit, the outer boundary of this ``hull'' has the same
distribution as the outer boundary of  the hull of
chordal $SLE_6$ going from the vertex
$A_3 \cap A_1$ to the vertex $A_3 \cap A_2$.
The identification with $SLE$ comes from the
conformal invariance assumption; Schramm determined
the value $\kappa = 6$ from a particular crossing
probability, but we now understand this in terms
of the locality property which scaling limits of
these boundary curves can be seen to satisfy.
Cardy's formula (and generalizations) were computed
for $SLE_\kappa$ in \cite{LSW1}.

Recently Smirnov \cite{Smirnov} made a major
breakthrough by proving
 conformal
invariance and Cardy's formula for the limit of critical
percolation in the triangular lattice.   As a corollary, the
identification of the limit with $SLE_6$ has become a theorem.
This has also
led to  rigorous proofs of a number of critical
exponents for the lattice model \cite{LSWonearm,S2,
SW}.  The basic strategy is to compute the exponent for
$SLE_6$ and to then to use Smirnov's result to relate
this exponent to lattice percolation.

It is  an open problem to show that critical percolation
on other planar lattices, e.g., bond percolation on the square
lattice, has the same limiting behavior.

\subsection{Loop-erased random walk}\vskip-5mm \hspace{5mm }

{\em Loop-erased random walk (LERW)} in a finite set $A \subset
\Z^2$ starting at $0 \in A$ is the measure
on self-avoiding paths obtained from starting a simple
random walk at the origin, stopping at the first time
that it leaves $A$, and erasing loops chronologically
from the path.  It can also be defined as a nonMarkov chain
which at each time $n$ chooses a new step using
probabilities weighted by the probability that simple
random walk starting at the new point avoids the path
up to that point (see, e.g., \cite{LKesten}).  It is also related
to uniform spanning trees; if one choose a spanning
tree uniformly among all spanning trees of $A$, considered
as a graph with appropriate boundary conditions, then the
distribution of the unique self-avoiding path from the origin to the
boundary is the same as LERW.  Wilson gave a beautiful
algorithm to generate uniform spanning trees using LERW
\cite{Wilson}.

One can hope to define a scaling limit
of planar LERW on a domain connecting an
interior point to a boundary point by taking LERW on finer
and finer grids and taking the limit.  There are a number
of reasons to believe that this limit is conformally invariant.
For example. the limit of simple random walk (Brownian motion)
is conformally invariant and the ordering of points used in
the loop-erasing procedure is not changed under conformal
maps.  Also, certain crossing probabilities for LERW can be
given by determinants of probabilities for simple random walk
(see \cite{Fomin}),
and hence these quantities are conformally invariant.
Kenyon \cite{Kenyon} used a conformal invariance argument
(using a determinant relation from a related domino
tiling model)
to prove that the growth exponent for LERW is $5/4$; roughly,
this says it takes about $r^{5/4}$ steps for a LERW to travel
distance $r$.

Schramm \cite{S1} showed that under the assumption of
conformal invariance, the scaling limit of  LERW must be
radial $SLE_2$.  He used conformal invariance and a natural
Markovian-type property of LERW to conclude that it must
be an $SLE_\kappa$, and then he used Kenyon's result to determine
$\kappa$. Recently, Schramm, Werner, and I \cite{LSWlerw} proved
that the scaling limit of loop-erased random walk is $SLE_{2}$.

There is another path obtained from the uniform spanning tree that
has been called the uniform spanning tree Peano curve. This
path, which lies on the dual lattice, encodes the entire tree
(not just the path from the origin to the boundary).  A similar,
although somewhat more involved, argument can be used to show
that this process converges to the space-filling curve $SLE_8$
\cite{LSWlerw}.

\subsection{Self-avoiding walk}\vskip-5mm \hspace{5mm }

A {\em self-avoiding walk (SAW)} in the lattice $\Z^2$ is a nearest neighbor
walk with no self-intersections.  The problem of the SAW is
to understand the  uniform measure on all such walks of
a given length (or
sometimes the measure that assigns weight $a^n$ to all walks
of length $n$).  It is still an open problem to prove there
is a limiting distribution; it is believed that such a limit in
conformally invariant (see \cite{LSWsaw} for precise statements).
However, if the conjectures hold there is only one possible
limit, $SLE_{8/3}$.

The conformal invariance property leads one to conclude
that the limit must be an $SLE_\kappa$ and $\kappa
\leq 4$ is needed in order to have a measure on simple
paths.
The property that $SLE_{8/3}$ has that is not held by
$SLE_\kappa$ for other $\kappa \leq 4$ is the {\em restriction
property}.  The restriction property is similar to, but not
the same, as the locality property.  Let $A$ be a compact
subset and $\Phi$ the transformation
 as in \S\ref{section 2}
Then \cite{LSWrest} if
 $\gamma[0,\infty)$ is   an $SLE_{8/3}$ path
from $0$ to $\infty$,   the distribution of
$\Phi \circ \gamma$ given the event
$\{\gamma[0,\infty) \cap A = \emptyset \}$ is
the same as (a time change of) $SLE_{8/3}$.
In fact, the probability  that
$\{\gamma[0,\infty) \cap A = \emptyset \}$
is $\Phi'(0)^{5/8}$.

If the scaling limit of
SAW has a conformally invariant
limit then one can show easily that the limit
satisfies the restriction property.  Hence, the only candidate
for the limit (assuming a conformally invariant
scaling limit) is $SLE_{8/3}$.   The conjectures
for critical exponents for SAW can be interpreted
in terms of rigorous properties of $SLE_{8/3}$
(see \cite{LSWsaw}).  For
example, the Hausdorff dimension of $SLE_{8/3}$
paths is $4/3$ \cite{B,LSWrest}; this gives strong
evidence that the limit of SAWs should give paths
of dimension $4/3$.  Monte Carlo simulations \cite{Kennedy}
support the conjecture that the limit of SAW is
$SLE_{8/3}$.

\medskip

\noindent{\bf Acknowledgment.}  Oded Schramm and Wendelin Werner
should be considered co-authors of this paper since this describes
joint work.  I thank both of them for an exciting collaboration.

\label{lastpage}


\baselineskip 4.5mm

\begin{thebibliography}{aa}
\bibitem{Bass} R. Bass, {\em Probabilistic Techniques in
Analysis}, Springer-Verlag, 1995.
\bibitem{B} V. Beffara, Hausdorff dimensions for $SLE_6$,
preprint.
\bibitem{BL2} K. Burdzy \& G. Lawler, Non-intersection exponents
for random walk and Brownian motion. Part II: Estimates and applications
to a random fractal, {\it Ann. Probab}, 18 (1990), 981--1009.
\bibitem{Cardy} J. Cardy, Critical percolation in finite geometries,
{\it J. Phys. A}, 25 (1992), L201--L206.
\bibitem{DK1} B. Duplantier \& K.-H. Kwon, Conformal invariance
and intersections of random walks, {\it Phys. Rev. Lett.}, 61,
2514--2517.
\bibitem{D3} B. Duplantier, Random walks and quantum gravity in
two dimensions, {\it Phys. Rev. Let.} , 81 (1998), 5489--5492.
\bibitem{Fomin} S. Fomin, Loop-erased walks and total
positivity, {\it Trans. Amer. Math. Soc.}, 353 (2001),
3563--3583.
\bibitem{Kennedy} T. Kennedy, Monte Carlo tests of SLE
predictions for the 2D self-avoiding walk, {\it Phys. Rev. Lett.} 88
(2002),
130601
\bibitem{Kenyon} R. Kenyon, The asymptotic determinant of
the discrete Laplacian, {\it Acta Math.}, 185 (2000), 239--286.
\bibitem{Lang} R. Langlands, P. Pouliot, \& Y. Saint-Aubin,
Conformal invariance in two-dimensional percolation,
{\it Bull. Amer. Math. Soc. (N.S.)}  90 (1994), 1--61.
\bibitem{Lcut} G. Lawler, Hausdorff dimension of cut points for
Brownian motion, {\it Electronic J. Probab.}, 1 (1996),
paper no. 2.
\bibitem{Lfront} G. Lawler, The dimension of the frontier of planar
Brownian motion, {\it Electronic Comm. Probab.}, 1 (1996), paper
no. 5.
\bibitem{LBuda} G. Lawler, Geometric and fractal properties
of Brownian motion and random walk paths in two and three
dimensions, in {\it Random Walks, Budapest 1998}, Bolyai
Mathematical Studies, 9 (1999), 210--258.
\bibitem{LKesten} G. Lawler, Loop-erased random walk,
{\it Perplexing Problems in Probability}, Birkh\"auser (1999),
197--217.
\bibitem{LVienna} G. Lawler, An introduction to the
stochastic Loewner evolution, preprint.
\bibitem{LP} G. Lawler \& E. Puckette, The intersection
exponent for simple random walk, {\it Combinatorics, Probab.,
and Computing}, 9 (2000), 441--464.
\bibitem{LSW1} G. Lawler, O. Schramm,  \& W. Werner,
Values of Brownian intersection exponents I:
 Half-plane exponents, {\it Acta. Math.}, 187 (2001),
 237--273.
\bibitem{LSW2} G. Lawler, O. Schramm, \& W. Werner,
Values of Brownian intersection exponents II:
Plane exponents, {\it Acta. Math.}, 187 (2001),
275--308.
\bibitem{LSW3} G. Lawler, O. Schramm,  \& W. Werner,
Values of Brownian intersection exponents III: Two-sided
exponents, {\it Ann. Inst. Henri Poincar\'e}, 38 (2002),
109--123.
\bibitem{LSWanal} G. Lawler, O. Schramm,  \& W. Werner,
Analyticity of  intersection exponents for planar Brownian
motion, {\it Acta. Math.}, to appear.
\bibitem{LSWonearm} G. Lawler, O. Schramm, \& W. Werner,
One arm exponent for critical 2D percolation, {\it Electronic J.
Probab.}, 7 (2002), paper no. 2.
\bibitem{LSWlerw} G. Lawler, O. Schramm, \& W. Werner,
Conformal invariance of planar loop-erased random walk
and uniform spanning trees, preprint.
\bibitem{LSWsaw} G. Lawler, O. Schramm, \& W. Werner,
On the scaling limit of planar self-avoiding walk, preprint.
\bibitem{LSWrest} G. Lawler, O. Schramm, \& W. Werner,
Conformal restriction properties: the chordal case, in
preparation.
\bibitem{LW1} G. Lawler \& W. Werner,  Intersection exponents
for planar Brownian motion, {\it Annals of Probab.}, 27
(1999), 1601--1642.
\bibitem{LW2} G. Lawler \& W. Werner, Universality
for conformally invariant intersection exponents, J.
European Math. Soc. 2 (2000), 291--328.
\bibitem{mandel} B. Mandelbrot, {\it The Fractal Geometry of Nature},
Freeman, 1982.
\bibitem{SR} S. Rohde \& O. Schramm, Basic properties of SLE,
preprint.
\bibitem{S1} O. Schramm,
Scaling limits of loop-erased random walks and uniform
spanning trees, {\it Israel J. Math}, 118 (2001), 221--288.
\bibitem{S2} O. Schramm, A percolation formula, {\it Electronic
Comm. Probab.}, 8 (2001), paper no. 12.
\bibitem{Smirnov} S. Smirnov, Critical percolation in the plane:
Conformal invariance, Cardy's formula, scaling limits, {\it
C. R. Acad. Sci. Paris. Sr. I Math.}, 333 (2001), 239--244.
\bibitem{SW} S. Smirnov \& W. Werner, Critical exponents
for two-dimensional percolation, {\it Math. Res. Lett.}, to appear.
\bibitem{W} W. Werner, Critical exponents, conformal invariance
and Brownian motion, {\it Proceedings of the 3rd Europ. Congress
Math.}, Prog. Math 202 (2001), 87--103.
\bibitem{Wilson} D. Wilson, Generating random spanning trees
more quickly than the cover time, {\it Proceedings of the
Twenty-eighth Annual ACM Symposium on the Theory
of Computing}, ACM (1996), 296--303.
\end{thebibliography}
\end{document}